%%%%%%%%%%%%%%%%%%%%%%%%%%%%%%%%%%%%%%%%%%%%%%%%%%%%%%%%%%%%%%%%%%%%%%%%%%%
%% Arapura, Donu
%% 
%% Higgs line bundles, Green-Lazarsfeld sets,and maps of K\"ahler manifolds
%%   to curves
%% 
%% Let $X$ be a compact K\"ahler manifold. The set $\cha(X)$ of 
%%   one-dimensional complex valued characters of the fundamental group of 
%%   $X$ forms an algebraic group. Consider the subset of $\cha(X)$ 
%%   consisting of those characters for which the corresponding local 
%%   system has nontrivial cohomology in a given degree $d$. This set is 
%%   shown to be a union of finitely many components that are translates of
%%   algebraic subgroups of $\cha(X)$. When the degree $d$ equals 1, it is 
%%   shown that some of these components are pullbacks of the character 
%%   varieties of curves under holomorphic maps. As a corollary, it is 
%%   shown that the number of equivalence classes (under a natural 
%%   equivalence relation) of holomorphic maps, with connected fibers, of 
%%   $X$ onto smooth curves of a fixed genus $>1$ is a topological 
%%   invariant of $X$. In fact it depends only on the fundamental group of 
%%   $X$.
%% 
%% publ:  Bull. Amer. Math. Soc. (N.S.) 26(1992) no. 2
%% pp:    310-314
%% type:  Research Announcement        markup: amstex    file size: 19K
%% contact:dvb@math.purdue.edu
%% 
%% copyright: American Math. Society copyright; see end of article
%% 
%% Include files necessary for this article: bull-ppt.tex
%% 
%%%%%%%%%%%%%%%%%%%%%%%%%%%%%%%%%%%%%%%%%%%%%%%%%%%%%%%%%%%%%%%%%%%%%%%%%%%
\input amstex 
\documentstyle{amsppt}
\input bull-ppt
\keyedby{bull283/rxa}
\define\Hom{\operatorname{Hom}}

\define\cha{\operatorname{char}}
\define\Higgs{\operatorname{Higgs}}
\define\Pic{\operatorname{Pic}}
\define\im{\operatorname{im}}
\define\tors{\operatorname{tors}}

\topmatter
\cvol{26}
\cvolyear{1992}
\cmonth{April}
\cyear{1992}
\cvolno{2}
\cpgs{310-314}
\title Higgs line bundles, Green-Lazarsfeld sets,\\and 
maps of K\"ahler
manifolds to curves\endtitle
\author Donu Arapura\endauthor
%\shortauthor{}
\shorttitle{Higgs line bundles and Green-Lazarsfeld sets}
\address Department of Mathematics, Purdue University, 
West Lafayette, Indiana
47907\endaddress
\ml dvb\@math.purdue.edu\endml
\date May 22, 1991\enddate
%\keywords{}
\abstract Let $X$ be a compact K\"ahler manifold. The set 
$\cha(X)$ of
one-dimensional complex valued characters of the 
fundamental group of $X$ forms
an algebraic group. Consider the subset of $\cha(X)$ 
consisting of those
characters for which the corresponding local system has 
nontrivial cohomology
in a given degree $d$. This set is shown to be a union of 
finitely many
components that are translates of algebraic subgroups of 
$\cha(X)$. When the
degree $d$ equals 1, it is shown that some of these 
components are pullbacks of
the character varieties of curves under holomorphic maps. 
As a corollary, it is
shown that the number of equivalence classes (under a 
natural equivalence
relation) of holomorphic maps, with connected fibers, of 
$X$ onto smooth curves
of a fixed genus $>1$ is a topological invariant of $X$. 
In fact it depends
only on the fundamental group of $X$.\endabstract
\thanks Partially supported by NSF\endthanks
\endtopmatter
\document
Let $X$ denote a compact K\"ahler manifold. Call two 
holomorphic maps 
$f\:X\to C$ and $f^{\prime}\:X\to C^{\prime}$, where $C$ 
and $C^{\prime}$ 
are curves, equivalent if there is an isomorphism 
$\sigma\:C\to C^{\prime}$ 
such that $f^{\prime}=\sigma\circ f$. Fix an integer 
$g>1$, and consider the 
set of equivalence classes of surjective holomorphic maps, 
with connected 
fibers, of $X$ onto smooth curves of genus $g$. We will 
see that this set is 
finite and that its cardinality $N_g(X)$ depends only on 
the fundamental 
group of $X$.

This result is deduced from a structure theorem for 
certain homologically 
defined sets of characters. A character of $X$ is a 
homomorphism of $\pi_1(X)$
into $\Bbb C^{\ast}$; it is unitary if the image of 
$\pi_1(X)$ lies in the 
unit circle $U(1)$. The set $\cha(X)$ of characters forms 
an affine 
algebraic group. For every character $\varrho\in\cha(X)$, 
we let 
$\Bbb C_{\varrho}$ denote the local system or locally 
constant sheaf on $X$ 
whose monodromy representation is given by $\varrho$. For 
each pair of 
integers $i$ and $m$, we define the subset
$\Sigma_m^i(X)$ of $\cha(X)$ to consist of those 
characters $\varrho$ for
which $\dim H^i(X,\Bbb C_{\varrho})\ge m$. We will denote 
$\Sigma_1^i(X)$ by $
\Sigma^i(X)$, and we will suppress the dependence on $X$ 
when there is no
danger of confusion. We will call a subset $S$ of 
$\cha(X)$ a unitary translate
of an affine subtorus if there exists a unitary character 
$\varrho\in\cha(X)$
such that $\varrho S$ is a connected algebraic subgroup.
\thm{Theorem 1} For $X$, $i$, and $m$ as above, the set 
$\Sigma_m^i$ is a union
of finitely many unitary translates of affine subtori.
\ethm
By a component of $\Sigma_m^i$, we will mean a unitary 
translate of an affine
subtorus $T\subseteq\Sigma_m^i$ that is maximal with 
respect to inclusion.
Using results of Beauville \cite{B1}, \cite{B2}, we can 
explicitly describe the
positive dimensional components of $\Sigma^1$.
\thm{Theorem 2} Any positive-dimensional component of 
$\Sigma^1$ is a
translate of an affine subtorus by a torsion element in 
$\cha(X)$. If
$T\subseteq\Sigma^1$ is a positive-dimensional component 
containing the trivial
character, then there exists a surjective holomorphic map 
with connected fibers
$f\:X\to C$ onto a smooth curve of genus at least two such 
that
$T=f^{\ast}\cha(C)$\ethm
\thm{Corollary} If $g\ge2$ then $N_g(X)$ is finite and it 
depends only
$\pi_1(X)$. In other words, if $X^{\prime}$ is another 
compact K\"ahler
manifold with $\pi_1(X^{\prime})\cong\pi_1(X)$ then
$N_g(X^{\prime})=N_g(X)$.
\ethm
\demo{Sketch of proof} Using the theorem, we see that 
$N_g(X)$ counts the
number of $2g$-dimensional components of $\Sigma^1(X)$ 
containing the trivial
character. $\Sigma^1$ has a purely group theoretic 
description:
$\varrho\in\Sigma^1(X)$ if and only if $H^1(\pi_1(X),\Bbb 
C_{\varrho})\neq0$.
Therefore, an isomorphism 
$\varphi\colon\pi_1(X)\cong\pi_1(X^{\prime})$ induces a
bijection $\varphi^{\ast}\:\cha(X^{\prime})\to\cha(X)$ 
such that
$\varphi^{\ast}(\Sigma^1(X'))=\Sigma^1(\roman X).\quad\qed$
\enddemo
Using Hodge theory, we can give a different, more analytic 
description of
$\cha(X)$. By a Higgs line bundle, we mean a pair 
$(L,\theta)$ consisting of a
holomorphic line bundle $L$ whose first Chern class 
$c_1(L)$ lies in the
torsion subgroup $H^2(X,\Bbb Z)_{\tors}$, together with a 
holomorphic $1$-form
$\theta$. The set of Higgs lines bundles $\Higgs(X)$ can 
be endowed with the
structure of a complex Lie group by identifying it with 
the product of the
Picard torus  $\Pic^0(X),H^2(X,\Bbb Z)_{\tors}$ and the 
vector space of
holomorphic $1$-forms. We define a map 
$\psi\:\cha(X)\to\Higgs(X)$ as follows:
$\psi(\varrho)=(L_{\varrho},\theta_{\varrho})$, where 
$L_{\varrho}$ is the
holomorphic bundle whose sheaf of sections is $\Bbb 
C_{\varrho}\otimes_{\Bbb C}
O_X$ and $\theta_{\varrho}$ is the (1, 0) part of 
$\log\|\varrho\|$ viewed as a
cohomology class under the isomorphism $H^1(X,\Bbb 
R)\cong\Hom(\pi_1(X),
\Bbb R)$.  Then $\psi$ is an isomorphism of topological 
groups (but not of
complex Lie groups). Simpson \cite{S} introduced the 
concept of a $\Higgs$
bundle of arbitrary rank on a K\"ahler manifold; however, 
the notion of Higgs
line bundle also occurs implicitly in the work of Green 
and Lazarsfeld
\cite{GL1}, \cite{GL2} and Beauville.

Before describing the image of $\Sigma_m^i$ under $\psi$, 
we need to define the
cohomology group of a Higgs line bundle $(L,\theta)$
$$H^{pq}(L,\theta)=\frac{\ker(H^q(X,\Omega_X^p\otimes
L)\buildrel{\wedge\theta}\over\to H^q(X,\Omega_X^{p+
1}\otimes L))}
{\im(H^q(X,\Omega_X^{p-1}\otimes 
L)\buildrel{\wedge\theta}\over\to
H^q(X,\Omega_X^p\otimes L))}.
$$

The next theorem follows by combining the results of Green 
and Lazarsfeld
\cite{GL1, 3.7} with those of Simpson \cite{S, 3.2}.
\thm{Theorem 3} For each $i$ there is an isomorphism 
$$H^i(X,\Bbb C_\varrho)\cong\bigoplus_{p+
q=i}H^{pq}(\psi(\varrho)).
$$
\ethm
We define the sets
$$\align
\sigma_m^{pq}&=\{(L,\theta)\in\Higgs(X)|\dim 
H^{pq}(L,\theta)\ge m\},\\
S_m^{pq}&=\{L\in\Pic^0(X)|\dim H^q(X,\Omega_X^p\otimes 
L)\ge m\}.
\endalign
$$
The set $S_m^{pq}$ was defined by Green and Lazarsfeld; it 
equals the
intersection of $\sigma_m^{pq}$ with $\Pic^0(X)\times\{0\}$.
\thm{Corollary} $\psi(\Sigma_m^i)=\bigcup_\mu\bigcap_{0\le 
k\le
i}\sigma_{\mu(k)}^{k,i-k}$, where $\mu$ runs over all 
partitions of $m$, i.e.,
functions $\mu\:\{0\cdots i\}\to \{0,1,2,\dots\}$ such 
that $\Sigma\mu(k)=m$.
\ethm
Let $\Bbb R^+$ denote the set of positive real numbers 
viewed as a group under
multiplication. A number $t\in\Bbb R^+$ acts on a Higgs 
line bundle by the rule
$t\ \ast\ (L,\theta)=(L,t\theta)$. We can transfer this 
action to $\cha(X)$ via $
\psi$, namely, $t\ \ast\ \varrho=\psi^{-1}(t\ {\ast}\ 
\psi(\varrho))$.
After choosing generators for $\pi_1(X)$, we can identify 
the connected components of
$\cha(X)$ with a product of $\Bbb C^{\ast}$'s. Under this 
identification the 
$\Bbb R^+$ action is described by 
$$t\ {\ast}\ (r_1e^{i\lambda_1},r_2e^{i\lambda_2},%
\dots)=(r_1e^{it\lambda_1},r_2e^{it
\lambda_2},\dots)$$
where $r_1,r_2,\dots\lambda_1,\dots\in\Bbb R$.

We can now indicate the idea of the proof of the first 
theorem. Using a Cech
complex, it is possible to write down equations for 
$\Sigma_m^i$, so we
conclude that this is an algebraic subset of $\cha(X)$. 
The corollary to
Theorem 3 shows that this set is stable under the $\Bbb R^+
$ action. The
theorem now follows from
\thm{Proposition} If $V\subseteq(\Bbb C^{\ast})^n$ is a 
closed irreducible
subvariety stable under the above $\Bbb R^+$ action, then 
$V$ is a unitary
translate of an affine subtorus.
\ethm
\demo{Sketch of proof} The Zariski closure of any orbit 
$\Bbb R^+\ \ast\ v$, with
$v\in(\Bbb C^{\ast})^n$, can be shown to be a unitary 
translate of an affine
subtorus. One then checks that for a sufficiently general 
point $v\in V$, the
orbit $\Bbb R^+\ast v$ is Zariski dense in $V.\qed$
\enddemo
As a corollary to Theorem 1, we obtain a new proof of a 
theorem of Green and
Lazarsfeld \cite{GL2} about the structure of $S_m^{pq}$. 
We say that a subset
$T$ of the Picard group $\Pic(X)$ is a translate of a 
complex subtorus if there
is an element $\tau\in\Pic(X)$ such that $\tau+T$ is a 
connected complex Lie
subgroup.
\thm{Corollary} There exist a finite number of translates 
of complex subtori
$T_i$ of $\Pic(X)$ and subspaces $V_i$ of the
space of holomorphic $1$-forms on
$X$ with $\dim T_i=\dim V_i$, such that $\sigma_m^{pq}$ is 
a union of
$T_i\times V_i$. In particular $S_m^{pq}$ is the union of 
those $T_i$ contained
in $\Pic^0(X)$.
\ethm
\demo{Sketch of proof} $\sigma_m^{pq}$ is an analytic 
subvariety of $\Higgs(X)$.
Choose an irreducible component $U$ of this set. Let $i=p+
q$ and for
$k\in\{0,\dots i\}$ define 
$$\mu(k)=\max\{n|U\subseteq\sigma_n^{k,i-k}\}.
$$
Then $U$ is an irreducible component of 
$\bigcap_i\sigma_{\mu(k)}^{k,i-k}$
that
is not contained in 
$\bigcap_i\sigma_{\mu^{\prime}(k)}^{k,i-k}$ for any other
partition $\mu^{\prime}$ of $M=\sum_j\mu(j)$. Thus $U$ is 
an irreducible
component of $\psi(\Sigma_M^i)$. By the theorem, it can be 
shown that any
irreducible component of $\psi(\Sigma^i_M)$
is the image under $\psi$ of a unitary translate
of an affine subtorus; such a set is of the form $T\times 
V$, where $T$ is a translate of a
complex subtorus of $\Pic(X)$ and $V$ is a subspace of 
1-forms of the same
dimension.\qed
\enddemo
We will call an unramified cover of $X$ with abelian 
Galois group an abelian
cover. The maximal abelian cover $X^{\roman{ab}}$ is 
obtained as the quotient of the
universal cover by the commutator subgroup 
$\pi_1(X)^{\prime}$. The Galois
group of $X^{\roman{ab}}$ over $X$ is precisely 
$H_1(X,\Bbb Z)$. The homology groups
$H_i(X^{\roman{ab}},\Bbb Z)$ are finitely generated as 
$\Bbb Z[H_1(X,\Bbb Z)]$-modules
although not necessarily as abelian groups. Our next 
theorem give partial
support to some conjectures of Beauville \cite{B2} and 
Catanese \cite{C} on the
structure of Green-Lazarsfeld sets.
\thm{Theorem 4} Fix an integer $N$. Suppose that 
$H^i(X^{\roman{ab}},\Bbb Z)$ is a
finitely generated abelian group for all $i<N$. Then 
\roster
\item"$\roman{(a)\phantom{\prime}}$" $\Sigma^i(X)$ 
consists of a finite set of torsion points of
$\cha(X)$ whenever $i<N$.
\item"$\roman{(a^{\prime})}$" $S_1^{pq}(X)$ consists of a 
finite set of torsion points
in $\Pic^0(X)$ whenever $p+q<N$.
\item"$\roman{(b)\phantom{\prime}}$" There is a finite 
sheeted abelian cover $X^{\prime}\to X$ such that
$\Sigma^i(X^{\prime})=\{1\}$ where $1$ is the trivial 
character whenever $i<N$.
\item"$\roman{(b^{\prime})}$" 
$S_1^{pq}(X^{\prime})=\{O_X\}$ whenever $p+q<N$.
\item"$\roman{(c)\phantom{\prime}}$" $\Sigma^N(X)$ has a 
positive-dimensional component if and only if
$H^N(X^{\roman{ab}},\Bbb Q)$ is infinite-dimensional.
\item"$\roman{(c^{\prime})}$" $S_1^{pq}(X)$ has a 
positive-dimensional component for
some $p$ and $q$, with $p+q=N$, if and only if 
$H^N(X^{\roman{ab}},\Bbb Q)$ is
infinite-dimensional.
\endroster
\ethm
\demo{Sketch of proof of $(\roman a)$} Let $V$ be a 
finite-dimensional $\Bbb C$-vector
space upon which $A=H_1(X,\Bbb Z)$ acts. A character 
$\varrho$ will be called a
weight of $V$ if there is a nonzero $v\in V$ such that for 
all $a\in
A,\,av=\varrho(a)v$. We prove a vanishing/nonvanishing 
theorem: $H^0(A,V\otimes
_{\Bbb C}\Bbb C_\varrho)=0$ if $\varrho^{-1}$ is a weight 
of $V$, otherwise
$H^p(A,V\otimes_{\Bbb C}\,\Bbb C_\varrho)=0$
for all $p$. Let $W$ be the union of the set of
weights of $H^i(X^{\roman{ab}},\Bbb 
C)=H^i(X^{\roman{ab}},\Bbb Z)\otimes\Bbb C$ with $i<N$,
and let $W^{-1}$ be the set of inverses of these weights. 
Associated to the
cover $X^{\roman{ab}}$ there is a spectral sequence
$$E_2^{pq}=H^p(A,H^q(X^{\roman{ab}},\Bbb C)\otimes_{\Bbb 
C}\Bbb C_\varrho)\Rightarrow H^{p+q}(X,
\Bbb C_{\varrho}).
$$
This together with the vanishing/nonvanishing theorem 
implies that
$\bigcup_{i<N}\Sigma^i(X)=W^{-1}$. Therefore the sets 
$\Sigma^i(X)$ are finite
when $i<N$, and so by Theorem 1 they must consist of 
unitary characters.

Let $K$ be the number field obtained by adjoining to $\Bbb 
Q$ all the
eigenvalues of generators of $A$ acting on 
$H^i(X^{\roman{ab}},\Bbb Z)$ with $i<N$.
Then $W$ is defined over the ring of integers $O_K$ of 
$K$. In other words
there is a subset 
$W^{\prime}\subset\Hom(\pi_1(X),O_K^{\ast})$ such that
$W=\bigcup_{i\colon K\to\Bbb C}i(W^{\prime})$. Since we 
have shown that the
characters in $W$ are also unitary, it follows by a 
theorem of Kronecker that
they must have finite order.$\quad$\qed
\enddemo
\thm{Corollary} The following are equivalent.
\roster
\item"(a)" $H_1(\pi_1(X)^{\prime},\Bbb Q)$ is 
infinite-dimensional.
\item"(b)" There is a finite sheeted abelian cover of $X$ 
that maps onto a
curve of genus at least two.
\endroster
\ethm
\demo{Sketch of proof of $(\roman a)\Rightarrow(\roman 
b)$} If $H_1(\pi_1(X)^{\prime},
\Bbb Q)\cong
H_1(X^{\operatorname{ab}},\Bbb Q)\cong 
H^1(X^{\operatorname{ab}},\Bbb Q)$ is infinite-dimensional 
then $\Sigma^1(X)$ has a
positive-dimensional component. By theorem 2, this 
component is a translate of
an affine subtorus by a torsion element. Therefore there 
is a finite abelian
cover $X^{\prime}$ of $X$ such that the pull back of this 
component, which lies
in $\Sigma^1(X^{\prime})$, contains the trivial character. 
Then Theorem 2 shows
that $X^{\prime}$ maps onto a curve of genus at least 
$2.\quad$\qed
\enddemo
\heading Acknowledgments\endheading
I would like to thank A. Beauville, P. Bressler, M. Green, 
R. Hain, R.
Lazarsfeld, M. Nori, M. Ramachandran, and C. Simpson for 
helpful conversations
and correspondence.

\enddocument